\documentclass[12pt]{amsart}
\usepackage[letterpaper]{geometry}
\usepackage[utf8]{inputenc}
\usepackage{graphicx}
\usepackage{url}
\usepackage[table]{xcolor}

\begin{document}

\title{Bibliometric Analysis of Senior US Mathematics Faculty}

\author{Joshua Paik}
\address{School of Mathematics and Statistics, University of St Andrews}
\curraddr{Department of Mathematics, The Pennsylvania State University}
\email{joshdpaik@gmail.com}

\author{Igor Rivin}
\address{Mathematics Department, Temple University}
\curraddr{Research Department, Edgestream Partners, LLP}
\email{rivin@temple.edu}
\keywords{citations, ranking}
\subjclass{00A99}
\begin{abstract}
We introduce a methodology to analyze citation metrics across fields of Mathematics. We use this methodology 
to collect and analyze the MathSciNet (\url{http://www.ams.org/mathscinet}) profiles of Full Professors of Mathematics at all 131 R1, research oriented US universities. The data recorded was citations, field, and time since first publication. We perform basic analysis and provide a ranking of US math departments, based on age corrected and field adjusted citations. 
\end{abstract}

\maketitle

\section*{Introduction}
For as long as the authors can remember, there has been discussion of comparable quality of various researchers (in all fields of research, but the authors are most familiar with mathematics, so this paper concerns itself \textbf{with mathematics exclusively}). While such a comparison is not strictly speaking possible (mathematics is not like competitive swimming, where a single number determines the better swimmer), those of us who have been on hiring committees have needed to compare researchers in diverse fields, and those of us who have had students (or job offers) have had to have some sort of estimate of the quality of people independently of age and field (and building on that, to have some reasonably gauge of the quality of \textit{departments}). The (admittedly ambitious) purpose of this note is to propose an \emph{objective} metric, based entirely on citations data (as such, it can be gamed, as can any metric be). 

Briefly, we normalize the number of MathSciNet citations by dividing it by the number of years since the author's first paper raised to the magical power $1.3$. We further segment mathematics into a number of ``major fields", assign mathematicians to fields (this is very difficult for some people, including, ironically, the authors of this paper), and compute the $z$-score of the normalized citation number. For each department we then compute the mean $z$-score of faculty to compute the department's ranking.
\section{Data}

Data was collected between January 16 - January 24, 2020. After accessing the list of R1 schools, we found the faculty lists from the relevant departmental web page, determined their level of seniority, and searched their profiles on MathSciNet. In total we collected citation records for 2807 math professors at 131 different institutions. We then collected the total citations, the year of earliest indexed publication, and the field of the most cited publication for each mathematician. A second pass through the data set occurred between January 30 - February 2, 2020, and discovered errors were corrected. This data set is as complete as possible. 

There were initially $65$ fields that were most cited as classified by MathSciNet, which we reduced to $20$ fields, using the mapping in the appendix. This mapping was constructed using general expertise on the way each field worked, and is recorded in the appendix. 

\section{Exploratory Data Analysis}
\subsection{Distribution}
Citations and Citations/Year$^{1.3}$ appear exponentially distributed after transformation by a square root. 
\begin{center}
\includegraphics[scale=0.4]{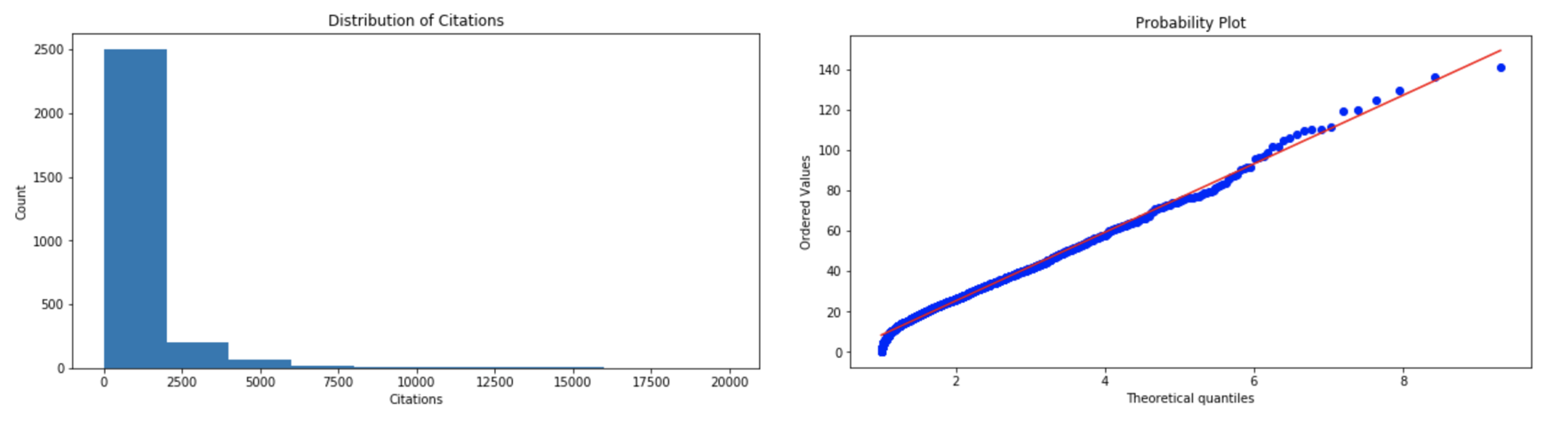}

Figure 1: Distribution of citations and an exponential qq-plot of square root citations.
\end{center}

\begin{center}
    \includegraphics[scale=0.4]{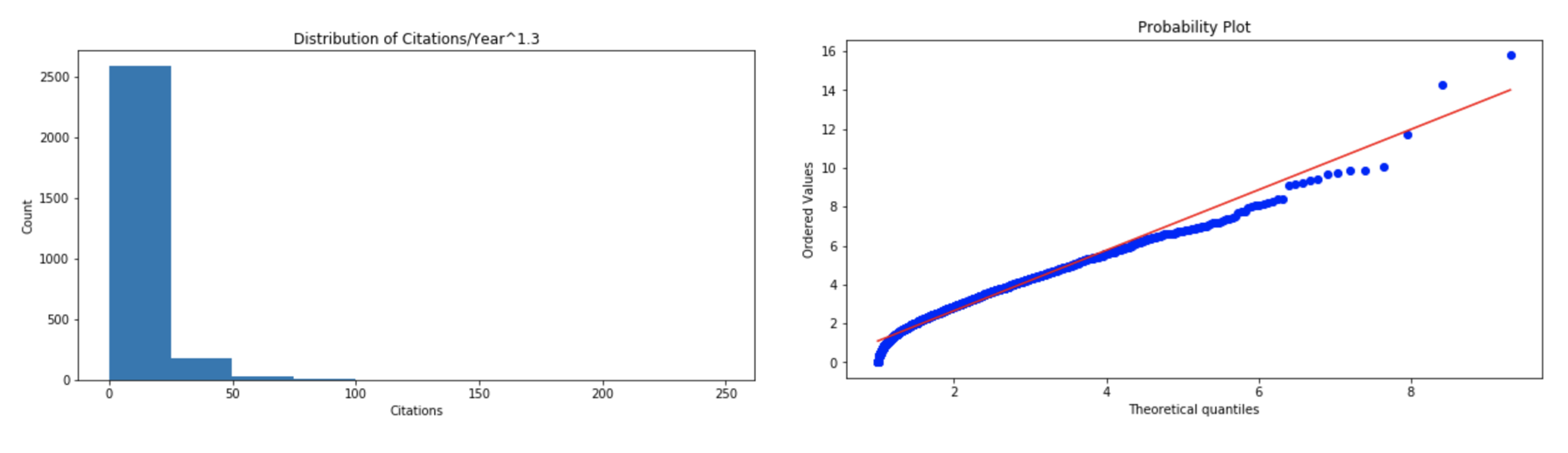}
\end{center}

Figure 2: Distribution of citations per year$^{1.3}$ and an exponential qq-plot of square root citations per year$^{1.3}$.

\subsection{Citations/Year$^{1.3}$ vs. Age}

As noted in earlier work \cite{paik2020data}, citations and year are positively correlated, but citations per year$^{1.3}$ and year are not correlated\footnote{The exponent $1.3$ was determined \it{via} a statistical analysis}. We repeat this analysis to check for robustness and find that when linearly regressed, citations/year$^{1.3}$ and year have a slope of $0.0338$ with a 95\% confidence interval of $[-0.00899, 0.0766]$. The $p$-value is $0.122$, so we fail to reject the null hypothesis that the slope is zero, and $R^2 = 0.03$. We conclude that citations/year$^{1.3}$ and year are not correlated.

\begin{center}
    \includegraphics[scale=0.5]{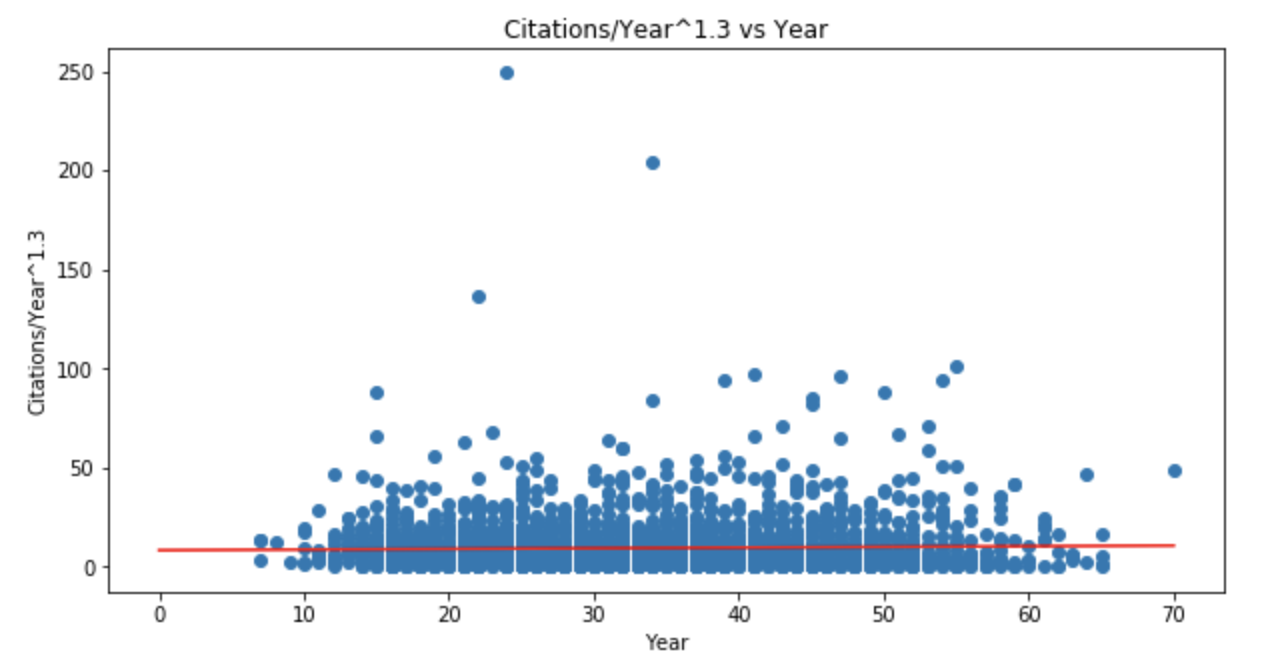}
    
Figure 3: Scatter plot of citations/year$^{1.3}$ and Years elapsed since first publication. The red line indicates the regressed line with equation $0.0338 Year + 8.40 = Citations/Year^{1.3}$. The $p$-value for slope is $0.122$ and the $R^2 = 0.03$. We fail to reject the null hypothesis that the slope is zero and conclude they are not correlated. 
\end{center}

\subsection{Factored Linear Models}

Before proceeding with the analysis, we should assess the importance of the  explanatory variables when looking at differences in citations between mathematicians. We do this by constructing nested linear models with all three combinations of \textit{Year} and \textit{Field}, and determine the best model with Akaike Information Criterion, AIC \cite{bishop2006machine}. Per standard interpretation the lower the AIC, the better the model. Let $C =$Citations, $a =$ Age, and $f = $ Field. 

\begin{center}
 \begin{tabular}{||c c ||} 
 \hline
 Model & AIC score \\ [0.5ex] 
 \hline\hline
 $\log(C) = \beta_0 a + \beta_1 f+ \epsilon$ & $9122.124$ \\
 \hline
 $\log(C) = \beta_0 f + \epsilon$ & $9347.318$  \\ 
  \hline
 $\log(C) = \beta_0 a + \epsilon$ & $9514.644$  \\ 
 \hline
\end{tabular}

Figure 4: Table of AIC scores for tested linear models. As the AIC score is lowest, the model consisting of both Age and Field is best, where Age impacts citations positively. For more detailed information on the model, refer to the ANOVA tables in the RMarkdown on Github. It is clear that certain fields contribute negatively to overall citations and other fields contribute positively. 
\end{center}

While it is not appropriate to pick a model for the sole reason that it minimizes AIC, it makes sense to consider both age and field. 

\subsection{Fields}

Different fields in mathematics have different citation practices. Some fields like Partial Differential Equations have more mathematicians, whereas some fields like Number Theory have fewer mathematicians. Some results from fields like topology are widely applicable across disciplines, whereas more obscure results are not. We quantify the bibliographic differences between fields. Note that the major fields below are larger categories containing potentially multiple MathSciNet tags, and the mappings are recorded in the appendix. 

\begin{center}
 \begin{tabular}{||c c c c c ||} 
 \hline
 Major Field & Mean Citations & S.D. Citations & Mean Cit/Year$^{1.3}$ & Count \\ [0.5ex] 
 \hline\hline
 PDE & 1472.07 & 2182.45 & 14.58 & 372\\ 
 \hline
 Computer Science & 1260.44 & 2223.06 & 14.08 & 225\\ 
 \hline
 Probability & 1165.92 & 1401.07 & 12.06 & 137\\ 
 \hline
 Harmonic Analysis & 1120.12 & 1336.62 & 10.51 & 200\\ 
 \hline
 Combinatorics & 1023.24 & 1673.53 &10.08 & 116\\ 
 \hline
  Algebra & 934.42 & 1310.59 &9.12 & 220\\ 
  \hline
   Algebraic Geometry & 846.62 &1308.80 & 9.51 & 169\\  
   \hline
    Geometry & 890.68 & 1486.72 & 8.87 & 311\\
   \hline
 Number Theory & 742.66 &920.31 & 7.38 & 159\\ 
   \hline
 Dynamics & 560.44 & 555.17 & 7.33 & 68\\
   \hline
 Mathematical Physics & 643.01 &716.41 & 7.25 & 96\\ 
 \hline
 Analysis & 977.18 & 1951.95 & 7.15 & 45\\ 
  \hline
 Applied Mathematics & 646.60 &  976.98 &6.87 & 299\\ 
   \hline
 Group Theory & 686.38 & 1151.64 & 6.74 & 81\\ 
  \hline
 Logic & 634.00 & 690.14 & 6.32 & 55\\ 
   \hline
 Complex Analysis & 612.86 & 725.15 & 6.17 & 115\\ 
   \hline
 Lie Groups & 512.02 & 590.59 & 4.78 & 43\\ 
   \hline
 Statistics & 220.73 & 331.15 & 3.10 & 83\\ 
   \hline
 History & 74.0 & 104.65 & 0.677 & 2 \\ 
   \hline
 Other & 5 & 6.61 & 0.074 & 11\\ 
 \hline
\end{tabular}

Figure 5: Mean citations and citations per year$^{1.3}$ including counts, split by field, from top to bottom ranked by mean citations per year$^{1.3}$ to account for age. 

\end{center}

We ran a permutation test between each field to verify the observed partial order. We report the inconclusive differences ($p$-value greater than $0.05$) between fields when comparing citations per year$^{1.3}$ in the appendix.

\subsection{Ranking of Departments}

The above figures shows that comparing mathematicians in two different fields is akin to comparing apples and oranges. The cleanest way to standardize this is to compute an interfield z-score of the citations per year$^{1.3}$, and hence associating a ``rank" to each mathematician. Then we computed the mean of the interfield $z$-scores of each full professor at the respective institution. We report here the top 20 schools and record the remaining schools in the appendix\footnote{Where by "Schools" we mean "Mathematics departments" - for Universities with separate Pure and Applied math departments, the ranking will be different if the departments were to be combined.}

\begin{enumerate}
  \item Princeton University
  \item Harvard University
  \item Stanford University
  \item University of Chicago
  \item Columbia University in the City of New York
  \item Massachussetts Institute of Technology
  \item University of California, Los Angeles
  \item University of Miami
  \item Yale University
  \item Brown University
  \item University of California, Berkeley
  \item New York University
  \item University of Oregon
  \item California Institute of Technology
  \item Duke University
  \item Stony Brook University
  \item Rutgers University-New Brunswick
  \item University of Virginia
  \item Texas A\&M University
  \item Northwestern University

\end{enumerate}

\section{Conclusion}
The rankings based on our normalized $z$-score (call it the PR score) correspond reasonably well with the ``folk" rankings of mathematicians. While we do not want to flatter or insult individuals by giving their scores here, we do give a ranking of departments, and we see that it, again, corresponds well with the ``folk" rankings. If they do not, we encourage the reader to look at the faculty pages of the departments in question. It seems, therefore, that there is, indeed, a fully quantitative way to produce \textbf{meaningful} rankings which work at least in a statistical sense - they fail for polymaths, and they also are less successful for mathematicians the bulk of whose work is not indexed by MathSciNet - in particular those who do interdisciplinary work.

\bibliographystyle{plain}
\bibliography{bibliography}

\section{Appendix}

\subsection{Code and Data}

Available at https://github.com/joshp112358/Differences .

\subsection{Classifications}
\begin{center}
    \begin{tabular}{||c | c||}
    \hline
         Major Field & Sub Fields  \\
        \hline\hline
         Algebra & Algebraic Topology; Associative \\
         & Rings and Algebra; Category theory, Homological algebra ;\\
         & Commutative rings and algebras; Field theory; \\
         & General algebraic systems; K-theory; Linear \\
         & and Multilinear Algebra, matrix theory; Associative \\
         & rings and algebras; Order, lattices,  ordered algebraic \\
         & structures. \\
    \hline
    Algebraic Geometry & Algebraic Geometry \\
    \hline
    Analysis & Difference and functional equations; \\
    & Integral equations; Integral transforms, operational calculus; \\
    &Ordinary differential equations; Real functions;\\
    &Special functions.\\
    \hline
    Applied Mathematics & Approximations and expansions; Biology \\
    & other natural sciences; Calculus of variations and optimal\\
    & control, optimization; Fluid mechanics,\\
    & Game theory, economics, social and behavioral sciences; \\
    & Geophysics, Mechanics of deformable sciences; \\
    & Mechanics of solids, Operations research, mathematical\\ 
    & programming; Systems theory, control.\\
    \hline
    Combinatorics & Combinatorics \\
    \hline
    Complex Analysis & Functions of a complex variable; Potential theory; Several \\
    & complex variables and analytic spaces\\
    \hline
    Computer Science & Computer Science; Numerical Analysis; \\
    & Information and communication, circuits. \\
    \hline
    Dynamics & Dynamical Systems and Ergodic Theory \\
    \hline
    Geometry & Convex and discrete geometry; Differential Geometry; \\
    & General topology; Geometry; \\
    & Manifolds and cell complexes;\\
    \hline
    Group theory & Group theory and generalizations. \\
    \hline 
    Harmonic analysis & Abstract harmonic analysis; Fourier analysis; \\
    & Functional analysis; Global analysis, analysis on manifolds; \\
    & Measure and integration, Operator theory. \\
    \hline
    History & History and biography. \\
    \hline
    Lie Groups & Topological Groups, Lie Groups. \\
    \hline
    Logic & Logic and foundations; Mathematical logic and foundations;\\ 
    & Set theory\\
    \hline 
    \end{tabular}
\end{center}
\begin{center}
    \begin{tabular}{||c | c||}
    \hline
         Major Field & Sub Fields  \\
        \hline\hline
            Mathematical Physics & Classical thermodynamics, heat transfer;\\
    & Mechanics of particles and systems; Optics, electromagnetic\\
    & theory; Quantum theory; Relativity and gravitational theory; \\
    & Statistical mechanics, structure of matter. \\
    \hline
    Number Theory & Number Theory \\
    \hline
    Other & Other \\
    \hline
        PDEs & Partial Differential Equations; Global Analysis,  \\
        & Analysis on manifolds\\
        \hline 
        Probability & Probability theory and stochastic processes \\
        \hline 
        Statistics & Statistics\\
        \hline
    \end{tabular}
\end{center}
\subsection{Ranking of Departments}

\begin{enumerate}
  \item Princeton University
  \item Harvard University
  \item Stanford University
  \item University of Chicago
  \item Columbia University in the City of New York
  \item Massachussetts Institute of Technology
  \item University of California, Los Angeles
  \item University of Miami
  \item Yale University
  \item Brown University
  \item University of California, Berkeley
  \item New York University
  \item University of Oregon
  \item California Institute of Technology
  \item Duke University
  \item Stony Brook University
  \item Rutgers University-New Brunswick
  \item University of Virginia
  \item Texas A\&M University
  \item Northwestern University
  \item University of Michigan
  \item Rice University
  \item The University of Texas at Austin
  \item Carnegie Mellon University
  \item University of Illinois at Chicago
  \item University of California, Irvine
  \item University of Pittsburgh
  \item Georgia Institute of Technology
  \item University of Minnesota
  \item Vanderbilt University
  \item Indiana University Bloomington
  \item SUNY at Albany
  \item University of California, San Diego
  \item University of North Texas
  \item University of Washington
  \item University of Connecticut
  \item Arizona State University
  \item Pennsylvania State University
  \item University of Southern California 
\item Purdue University                             
\item University of Illinois at Urbana-Champaign    
\item Cornell University                            
\item University of Maryland - College Park         
\item University of Utah                            
\item North Carolina State University               
\item Johns Hopkins University                      
\item University of California, Riverside           
\item University of California, Santa Cruz         
\item Washington University in St. Louis            
\item Wayne State University    
\item University of Pennsylvania                                 
\item Brandeis University                                        
\item Colorado State University                                  
\item University of Notre Dame                                   
\item University of California, Santa Barbara                    
\item University of North Carolina at Chapel Hill                
\item University of Houston                                      
\item University of Iowa                                         
\item The Ohio State University                                 
\item University of South Florida                               
\item Michigan State University                                 
\item University of California, Davis                          
\item Virginia Polytechnic Institute and State University       
\item University of Missouri                                    
\item University of Wisconsin - Madison                         
\item University of Massachusetts Amherst                       
\item University of South Carolina                              
\item Emory University                                          
\item University of Central Florida                             
\item University of Kentucky                                    
\item University of Florida                                     
\item University of Delaware                                    
\item Louisiana State University 
\item Syracuse University                                       
\item Georgia State University                                  
\item University of Colorado Denver                             
\item Boston University                                         
\item Tulane University of Louisiana                            
\item Clemson University                                        
\item University of Kansas                                      
\item University of Southern Mississippi                        
\item Boston College                                            
\item Mississippi State University                              
\item University of Rochester                                   
\item CUNY Graduate School and University Center                
\item The University of Tennessee                               
\item George Washington University                              
\item Georgetown University                                     
\item Florida State Universty                                   
\item Iowa State University                                     
\item University at Buffalo                                     
\item Northeastern University                                   
\item Tufts University                                          
\item University of Nebraska-Lincoln                            
\item University of Georgia                                     
\item University of New Hampshire                               
\item Virgina Commonwealth                                      
\item University of Cincinatti                                  
\item Dartmouth College                                     
\item Rennselaer Polytechnic Institute                      
\item University of Nevada, Reno            
\item West Virginia University           
\item Auburn University                    
\item The University of Texas at Arlington  
\item Texas Tech University                
\item University of Arizona                 
\item Binghamton University               
\item University of New Mexico              
\item The University of Alabama           
\item The University of Texas at Dallas    
\item George Mason University              
\item Florida Institute University         
\item University of Oklahoma                
\item University of Colorado Boulder        
\item University of Hawaii at Manoa         
\item Case Western Reserve University       
\item University of Alabama at Birmingham   
\item Oklahoma State University             
\item Kansas State University              
\item Temple University                    
\item Oregon State University               
\item Drexel University                   
\item University of Louisville             
\item University of Nevada, Las Vegas       
\item University of Wisconsin - Milwaukee   
\item Washington State University          
\item New Jersey Institute of Technology    
\item The University of Texas at El Paso   
\item University of Mississippi            
\item University of Arkansas              
\item Montana State University            
\end{enumerate}

\subsection{Inconclusive Permutation Tests between Fields}

We report the results of a one sided permutation test, when comparing cit/year$^{1.3}$ which failed to be significant at the $0.05$ level. We record the hypothesis on the left and the $p$-value to the right.

PDE $\geq$ Computer Science --- $p$-value: 0.397

PDE $\geq$ Probability --- $p$-value: 0.0768

Computer Science $\geq$ Probability --- $p$-value: 0.156

Probability $\geq$ Harmonic Analysis --- $p$-value: 0.113

Probability $\geq$ Combinatorics --- $p$-value: 0.1049

Harmonic Analysis $\geq$ Combinatorics --- $p$-value: 0.3824

Harmonic Analysis $\geq$ Algebra --- $p$-value: 0.0961

Harmonic Analysis $\geq$ Algebraic Geometry --- $p$-value: 0.1807

Combinatorics $\geq$ Algebra --- $p$-value: 0.2181

Combinatorics $\geq$ Algebraic Geometry --- $p$-value: 0.3265

Combinatorics $\geq$ Geometry --- $p$-value: 0.1544

Combinatorics $\geq$ Analysis --- $p$-value: 0.0719

Algebra $\geq$ Algebraic Geometry --- $p$-value: 0.6461

Algebra $\geq$ Geometry --- $p$-value: 0.3989

Algebra $\geq$ Dynamics --- $p$-value: 0.0813

Algebra $\geq$ Mathematical Physics --- $p$-value: 0.0546

Algebra $\geq$ Analysis --- $p$-value: 0.1232

Algebraic Geometry $\geq$ Geometry --- $p$-value: 0.2533

Algebraic Geometry $\geq$ Analysis --- $p$-value: 0.0782

Geometry $\geq$ Number Theory --- $p$-value: 0.0534

Geometry $\geq$ Dynamics --- $p$-value: 0.1076

Geometry $\geq$ Mathematical Physics --- $p$-value: 0.0704

Geometry $\geq$ Analysis --- $p$-value: 0.1497

Number Theory $\geq$ Dynamics --- $p$-value: 0.4906

Number Theory $\geq$ Mathematical Physics --- $p$-value: 0.4514

Number Theory $\geq$ Analysis --- $p$-value: 0.4417

Number Theory $\geq$ Applied Mathematics --- $p$-value: 0.2717

Number Theory $\geq$ Group Theory --- $p$-value: 0.2637

Number Theory $\geq$ Logic --- $p$-value: 0.1557

Number Theory $\geq$ Complex Analysis --- $p$-value: 0.0717

Dynamics $\geq$ Mathematical Physics --- $p$-value: 0.4614

Dynamics $\geq$ Analysis --- $p$-value: 0.4609

Dynamics $\geq$ Applied Mathematics --- $p$-value: 0.3377

Dynamics $\geq$ Group Theory --- $p$-value: 0.2999

Dynamics $\geq$ Logic --- $p$-value: 0.1738

Dynamics $\geq$ Complex Analysis --- $p$-value: 0.1272

Mathematical Physics $\geq$ Analysis --- $p$-value: 0.4916

Mathematical Physics $\geq$ Applied Mathematics --- $p$-value: 0.3506

Mathematical Physics $\geq$ Group Theory --- $p$-value: 0.337

Mathematical Physics $\geq$ Logic --- $p$-value: 0.23

Mathematical Physics $\geq$ Complex Analysis --- $p$-value: 0.1423

Mathematical Physics $\geq$ History --- $p$-value: 0.0525

Analysis $\geq$ Applied Mathematics --- $p$-value: 0.4035

Analysis $\geq$ Group Theory --- $p$-value: 0.3972

Analysis $\geq$ Logic --- $p$-value: 0.3154

Analysis $\geq$ Complex Analysis --- $p$-value: 0.2375

Analysis $\geq$ Lie Groups --- $p$-value: 0.1041

Analysis $\geq$ History --- $p$-value: 0.0881

Applied Mathematics $\geq$ Group Theory --- $p$-value: 0.4655

Applied Mathematics $\geq$ Logic --- $p$-value: 0.3609

Applied Mathematics $\geq$ Complex Analysis --- $p$-value: 0.2362

Applied Mathematics $\geq$ Lie Groups --- $p$-value: 0.0616

Applied Mathematics $\geq$ History --- $p$-value: 0.0534

Group Theory $\geq$ Logic --- $p$-value: 0.3728

Group Theory $\geq$ Complex Analysis --- $p$-value: 0.2834

Group Theory $\geq$ Lie Groups --- $p$-value: 0.0513

Logic $\geq$ Complex Analysis --- $p$-value: 0.4276

Logic $\geq$ Lie Groups --- $p$-value: 0.0613

Complex Analysis $\geq$ Lie Groups --- $p$-value: 0.0871

Lie Groups $\geq$ History --- $p$-value: 0.0531

Statistics $\geq$ History --- $p$-value: 0.2517

History $\geq$ Other --- $p$-value: 0.1533

\end{document}